 \documentstyle[12pt]{article}
\title{\bf Braided  m-Lie Algebras}
\author{
Shouchuan Zhang $^{a,~b}$ and   Yao-Zhong Zhang $^b$   \ \  \ \\
$a$. Department  of Mathematics,
Hunan University\\ Changsha  410082, \ P.R. China \\
$b$. Department of Mathematics, University of Queensland\\
Brisbane 4072, Australia\\ }
\date{}
\begin{document}
\newtheorem{Theorem}{\quad Theorem}[section]
\newtheorem{Proposition}[Theorem]{\quad Proposition}
\newtheorem{Definition}[Theorem]{\quad Definition}
\newtheorem{Corollary}[Theorem]{\quad Corollary}
\newtheorem{Lemma}[Theorem]{\quad Lemma}
\newtheorem{Example}[Theorem]{\quad Example}
\maketitle \addtocounter{section}{-1}
\begin {abstract}
Braided m-Lie algebras induced by multiplication are introduced, which
generalize Lie algebras, Lie color algebras and quantum Lie algebras.
The necessary and sufficient conditions for the braided m-Lie  algebras to
be strict Jacobi braided Lie algebras are given. Two classes of braided
m-Lie algebras are given, which are  generalized matrix braided m-Lie
algebras and braided m-Lie subalgebras of $End _F M$, where $M$ is
a Yetter-Drinfeld module over $B$ with dim $B< \infty $ .
In particular, generalized classical braided m-Lie algebras
$sl_{q, f}( GM_G(A),  F)$  and  $osp_{q, t} (GM_G(A), M, F)$  of generalized
matrix algebra $GM_G(A)$ are constructed  and
their connection with  special generalized matrix Lie superalgebra
$sl_{s, f}( GM_{{\bf Z}_2}(A^s),  F)$  and
orthosymplectic generalized matrix Lie super algebra
$osp_{s, t} (GM_{{\bf Z}_2}(A^s), M^s, F)$  are established.
The relationship between representations of braided m-Lie algebras and
their  associated algebras are established.

\end {abstract}

\section {Introduction}
The theory of Lie superalgebras has been developed systematically,
which includes  the representation theory and classifications of
simple Lie superalgebras and their varieties \cite {Ka77} \cite
{BMZP92}.  In many physical applications or in pure mathematical
interest, one has to consider not only ${\bf Z}_2$- or ${\bf Z}$-
grading but also $G$-grading of Lie algebras, where $G$ is an
abelian group equipped with a skew symmetric bilinear form given
by a 2-cocycle. Lie algebras in symmetric and more general
categories were discussed in \cite {Gu86} and \cite {GRR95}.
A sophisticated multilinear version of the Lie bracket was
considered in \cite {Kh99} \cite {Pa98}. Various generalized Lie
algebras have already appeared  under different names, e.g. Lie
color algebras, $\epsilon $ Lie algebras \cite {Sc79}, quantum and
braided Lie algebras,  generalized Lie algebras \cite {BFM96} and
$H$-Lie algebras \cite {BFM01}.

In  \cite {Ma94b}, Majid introduced braided Lie algebras from geometrical
point of view, which have attracted attention in
mathematics and mathematical physics (see e.g. \cite {Ma95} and references
therein).

In   this paper we introduce  braided m-Lie  algebras and (strict) Jacobi
braided Lie algebras from the algebraic point of view, rather than geometrical
view point as in \cite{Ma94b}. Our braided m-Lie algebras are different
from the braided Lie algebras defined by Majid \cite{Ma94b}. This is
verified by giving an example of braided m-Lie algebras
which is neither a braided Lie algebra (in the sense of
Majid) nor an Jacobi braided Lie algebra.  We give the necessary and
sufficient conditions for the braided m-Lie  algebras to be strict Jacobi
braided Lie algebras.  In section 2 we give $G$-gradings of
generalized matrix algebras
and construct braided m-Lie algebras  corresponding to such algebras,
which are called generalized matrix braided m-Lie algebras. This leads to the
braided m-Lie algebra
construction on path algebras, full matrix algebras and simple algebras.
 We construct  generalized classical braided m-Lie algebras of
generalized matrix algebras. In particular,
 special generalized matrix Lie color algebra $sl_{q, f}( GM_G(A),  F)$    and
ortho-symplectic generalized matrix Lie color algebra $osp_{q, t}
(GM_G(A), M, F)$  are related to the corresponding  Lie
superalgebras. In Section 3 we  give another class of braided
m-Lie algebras, i.e. braided m-Lie subalgebras of  $(End _FM)^-$.
We  show that representations of an  algebra $A$ associated to a
braided m-Lie algebra $L$ are also representations  of $L$.
Furthermore, we show that if $(M, \psi)$ is a faithful
representation of $L$, then representations of $End _FM$  are also
ones of  $L$.

Throughout, $ F$ is a field, $G$ is an additive group and
$r$ is a bicharacter of $G$;
$\mid $$ x $$\mid $ denotes the degree of $x$ and  $({\cal C}, C)$ is
a braided tensor category with braiding $C$. We  write $W \otimes f $ for
$id _W \otimes f$ and $f \otimes W$ for  $f \otimes id _W$.
Algebras discussed here may not have unity element.

\section {Braided  m-Lie Algebras}
In this section we introduce  braided m-Lie  algebras and (strict) Jacobi
braided Lie algebras.  We give the necessary and
sufficient conditions for the braided m-Lie  algebras to be strict Jacobi
braided Lie algebras.

\begin {Definition} \label {1.1}
Let  $(L, [\ \ ])$ be an object in the braided tensor category $
({\cal C }, C)$ with morphism $[\ \ ] : L \otimes L \rightarrow
L$. If there exists an algebra $(A, m)$ in  $ ({\cal C }, C)$ and
monomorphism $\phi : L \rightarrow A$ such that   $\phi [\  \ ] =
m (\phi \otimes \phi ) - m (\phi \otimes \phi )  C_{L, L},$   then
$(L, [ \ \  ])$ is called a braided m-Lie algebra in $ ({\cal C },
C)$ induced by multiplication of $A$ through $\phi $. Algebra $(A,
m) $ is called an algebra  associated to $(L, [ \ \ ])$.

\end {Definition}
A Lie  algebra is a braided m-Lie algebra in the category of ordinary vector
spaces, a Lie color algebra is a braided  m-Lie algebra
 in symmetric braided tensor category $ ({\cal M} ^{FG}, C^r)$ since the
 canonical map $\sigma: L \rightarrow U(L)$ is injective (see
\cite [Proposition 4.1]{Sc79}), a  quantum Lie  algebra is a braided
m-Lie algebra in the Yetter-Drinfeld category $ (^B_B{\cal YD}, C)$ by
\cite [Definition 2.1 and Lemma 2.2]{GM03}), and a ``good"
braided Lie  algebra is a braided  m-Lie algebra
 in the Yetter-Drinfeld category $ (^B_B{\cal YD}, C)$ by
 \cite [Definition 3.6 and Lemma 3.7]{GM03}).
For  a cotriangular Hopf algebra
$(H, r)$, the $(H,r)$-Lie algebra defined in \cite [4.1] {BFM01} is a
 braided m-Lie  algebra in the braided
tensor category $({}^H{\cal M}, C^r)$.  Therefore, the
braided m-Lie algebras generalize most known generalized Lie algebras.

For an algebra $(A, m)$ in $({\cal C}, C)$, obviously $L = A$ is  a braided
m-Lie algebra under operation $ [\ \ ] = m  - m   C_{L, L}$, which  is induced
by $A$ through $id _A$. This  braided  m-Lie algebra is written as $A^-$.

If $V$ is an object in ${\cal C}$ and $C_{V,V} = C_{V, V}^{-1},$ then we say
that the braiding is symmetric on $V$.

\begin {Example} \label {1.2} If $H$ is a  braided Hopf algebra in the
Yetter-Drinfeld module category $(^B_B {\cal YD}, C)$
with $B = FG$ and $C(x\otimes y) = r(\mid $$ y $$\mid , \mid $$ x $
$\mid )y \otimes x $ for any homogeneous elements $x, y \in H$,
then $P(H)=: \{ x \in H \mid x
\hbox { is a primitive element }\} $ is a braided m-Lie  algebra iff
the braiding $C$ is symmetric on $P(H)$.
\end {Example}

Indeed, it is easy to check  that $P(H)$ is the Yetter-Drinfeld
module. By simple computation we have $\Delta ([x, y ]) = [x, y]
\otimes 1 + 1 \otimes [x, y] + (1 - r(\mid $$x$$ \mid , \mid $$y$$
\mid )r(\mid $$y$$ \mid , \mid $$x$ $ \mid )) x \otimes y$ for any
homogeneous elements $x, y \in P(H).$
 Thus $[x, y] \in P(H)$ iff $r(\mid $$x$$ \mid , \mid $$y$$ \mid )
 r(\mid $$y$$ \mid , \mid $$x$$ \mid )=1,$ as asserted.

\begin {Theorem} \label {1.3} Let  $(L, [\ \ ])$ be a braided  m-Lie
 algebra in $({\cal C}, C)$.

 (i)  $(L, [\  \ ])$ satisfies the braided anti-symmetry (or quantum
 anti-symmetry):

 (BAS): $[\  \ ] = - [\  \ ] C_{L, L}$

\noindent  if and only if
 $mC_{L,L}= mC_{L,L}^{-1}.$

 (ii) If the braided anti-symmetry holds, then braided m-Lie  algebra
 $(L, [\  \  ], m)$ satisfies the (left) braided primitive  Jacobi identity:

  (BJ):\ \ $[\ \ ] (L\otimes [\ \ ]
)+[\ \ ] (L\otimes [\ \ ] ) (L  \otimes C_{L,L}^{-1})
(C_{L,L} \otimes L) + [\ \ ] (L \otimes [\ \ ]
)(C_{L,L}^{-1} \otimes L)(L \otimes C_{L,L} ))=0$, \\
 and the right  braided primitive Jacobi identity:

 (BJI'):\ \ $[\ \ ] ([\  \ ]\otimes L
)+[\ \ ] ([\  \ ]\otimes L ) (L \otimes C_{L,L}{}^{-1})
(C_{L,L} \otimes L) + [\ \ ] ([\  \ ]\otimes L)(C_{L,L}{}^{-1} \otimes L)
(L \otimes C_{L,L} ))=0.$
 \end {Theorem}

{\bf Proof .}
(i) Assume that $(L, [\  \ ])$ satisfies the braided anti-symmetry.
Since $[\  \ ] = - [\  \ ]C$ we have $m- m C = mCC - mC$ and $m =m CC,$
which implies $mC= m C^{-1}.$ The necessity is clear.

(ii) \begin {eqnarray*}
 \hbox { l.h.s. of  (BJI)} &=& m(L\otimes m)- m(L \otimes m)(L \otimes C) \\
& &- m C(L \otimes m) + m C(L \otimes m) (L \otimes C)\\
 & &+ m (L \otimes m) (L\otimes C^{-1}) (C \otimes L) \\
& &+ m (L \otimes m)(L \otimes C)(L \otimes C^{-1})(C \otimes L)\\
& &- m C (L \otimes m)(L \otimes C^{-1})(C \otimes L) \\
& &+ m C(L \otimes m)(L \otimes C)(L \otimes C^{-1})(C \otimes L)\\
& &+ m (L \otimes m) (C^{-1} \otimes L )(L \otimes C)\\
& &- m(L \otimes m)(L \otimes C)(C^{-1} \otimes L)(L \otimes C)\\
& &- mC(L \otimes m) (C^{-1} \otimes L)(L \otimes C)\\
& &+ m C(L \otimes m)(L \otimes C)(C^{-1} \otimes L)(L \otimes C) .
\end {eqnarray*}

We first check that $-$(6th term) =  12th term and 4th term = $-$(10th term).
 Indeed,\begin {eqnarray*}
\hbox { 12th term } &=& m C^{-1}(L \otimes m)
(L \otimes C^{-1})(C^{-1} \otimes L)(L \otimes C)\\
&=&m(m \otimes L)(L \otimes C^{-1})(C^{-1} \otimes L)(L \otimes C^{-1})
(C^{-1} \otimes L)(L \otimes C)\\
&=&m (L \otimes m)(L \otimes C)(C^{-1} \otimes L)
(L \otimes C^{-1})(C^{-1} \otimes L)(L \otimes C)\\
&=&m (L \otimes m)(C^{-1} \otimes L)(L \otimes C^{-1})
(C \otimes L)(C^{-1} \otimes L)(L \otimes C)\\
&=& - \hbox {(6th term )}.
\end {eqnarray*}
\begin {eqnarray*}
\hbox { 4th term }&=& m (m \otimes L)(C^{-1} \otimes L)
(L \otimes C^{-1})(C^{-1} \otimes L)\\
&=&m(L \otimes m)(L \otimes C^{-1})(C^{-1} \otimes L)(L \otimes C) \\
&=& - \hbox {(10th term )}.
\end {eqnarray*}
We can similarly show that 1st term = $-$ (11th term),
$-$(2nd term) = 8th term, $-$ (3rd term) = 5th term,
 $-$(7th term) = 9th term.
Consequently, $(BJI)$ holds. We can similarly show that $(BJI)'$ holds.  $\Box$

Readers can prove the above with the help of   braiding diagrams.

\begin {Definition} \label {1.4}  Let $[\  \ ]$  be a morphism from
$L\otimes L$  to $L$ in ${\cal C}$.
If $(BJI)$ holds, then  $(L, [\  \ ])$  is called a  (left )  Jacobi braided
Lie  algebra. If both of (BAS) and (BJI) hold
then $(L, [\  \ ])$  is called a  (left ) strict Jacobi braided  Lie  algebra.
\end {Definition}
Dually, we can define  right Jacobi braided  Lie  algebras and right strict
Jacobi braided Lie  algebras.
Left (strict) Jacobi braided Lie  algebras  are called  (strict) Jacobi
braided   Lie  algebras in general.

By Theorem \ref {1.3} we have
 \begin {Corollary} \label  {1.5} If  $(L , [\  \ ] )$ is a braided m-Lie
algebra, then the following conditions are equivalent:

(i)  $(L, [\  \ ])$  is a left (right) strict Jacobi braided Lie algebra.

(ii) $m C_{L,L} = m C_{L,L}^{-1}.$

(iii)  $[\  \ ] C_{L,L} = [\ \ ] C_{L,L}^{-1}.$

 \end {Corollary}

Furthermore, if $L$ is a space graded by $G$ with bicharacter $r$, then
the braided primitive Jacobi identity (BJI) becomes:
$$r (\mid c \mid , \mid a \mid) [a, [b, c]] +
 r (\mid b \mid , \mid a \mid) [b, [c, a]]+ r (\mid c \mid , \mid b \mid)
 [c, [a, b]] =0   \ \ \ \ \ \ \ \  \ \ \ \ \ \ \ \ (*) $$
for any homogeneous elements $a, b, c \in L.$  That is, $L$ is a Jacobi
braided Lie  algebra if and only if (*) holds.
For convenience, we let $J(a, b, c)$ denote the left hand side of $(*)$.

We now recall the (left ) braided Lie algebra defined by Majid
\cite [Definition 4.1]{Ma94b}.  A (left) braided Lie algebra
  in ${\cal C}$ is a coalgebra $(L, \Delta , \epsilon  )$ in ${\cal C}$,
 equipped with a morphism  $[\  \ ]: L \otimes L \rightarrow L$
 satisfying the axioms:

(L1)
$([\  \ ])(L \otimes [\  \ ])= [\  \ ]([\  \ ] \otimes
[\  \ ])(L \otimes C \otimes L)(\Delta \otimes L \otimes L)  $

(L2)  $C ([\  \ ] \otimes L) (L \otimes C)
(\Delta \otimes L) = (L \otimes [\  \ ])(\Delta \otimes L)$

(L3) $[\  \ ]$ is a coalgebra morphism in $({\cal C}, C).$

\noindent Axiom (L1) is called the left braided Jacobi identity.

There exist braided
m-Lie  algebras which are neither  braided Lie algebras nor Jacobi braided
Lie  algebras as is seen from the following example.
\begin {Example} \label {1.6} (see \cite  {Ma95} )
Let $L = F\{x\} / <x^n>$ be an algebra in $(^{F {\bf Z}_n} {\cal M}, C^r)$
with a primitive nth root $q$ of 1 and
$r (k, m) = q ^{km}$ for any $k, m \in {\bf Z}_n,$ where $n$ is a natural
number.

(i) If  $n > 3$,  then  $(L , [\  \ ])$ is  a braided m-Lie algebra but
is neither a Jacobi braided Lie  algebra nor a braided Lie algebra of Majid.

(ii) If  $n = 3$,  then  $(L , [\  \ ])$ is a braided m-Lie algebra and a
Jacobi braided Lie  algebra  but is neither a strict Jacobi braided Lie
algebra nor a braided Lie algebra of Majid.
\end {Example}

{\bf Proof.}
(i) Since $J(x, x, x) =3qx^3 ( 1 -q - q^2 + q^3) \not=0 $ we have that
$(L, [\  \ ] )$ is not a Jacobi braided Lie algebra.
If $(L, [\  \ ], \Delta , \epsilon )$ is a braided Lie algebra,
since $[\  \  ]$ is a coalgebra homomorphism, we have  that
\begin {eqnarray*}\label {e1}
\hbox {[1  1]} &=& 0, \hbox { implying } \epsilon (1)=0\\
\hbox {[ } x , x\hbox{ ]} &=& x^2 (1 - q),  \hbox { implying }
\epsilon (x^2) (1-q)= \epsilon (x)^2\\
\hbox {[}x , x^2\hbox {]} &=& x^3 (1 - q^2),  \hbox { implying }
\epsilon (x^3) (1-q^2)= \epsilon (x) \epsilon (x^2)\\
&\cdots & \\
\hbox {[}x , x^{n-1}\hbox {]} &=& x^n (1 - q^{n-1}),  \hbox { implying }
0 = \epsilon (x^n) (1-q^{n-1})= \epsilon (x) \epsilon (x^{n-1}).\\
\end {eqnarray*}
Thus $\epsilon (x^m ) = 0$ for $m =0, 1, 2, \cdots , n-1$, which contradicts
the fact that $(L, \Delta , \epsilon )$ is a coalgebra.

(ii) Since  $J(x^k,x^l,x^m)=0$ for  any $k, l , m \in {\bf Z}_3$, we
have that   $(L , [\  \ ])$ is a Jacobi braided Lie  algebra. It follows from
$[x, x] = x^2(1-q) \not= - q [x, x]$ that  $(L , [\  \ ])$ is not a strict
Jacobi braided Lie  algebra.
$\Box$

Note that  $L$ in the above example may never be a braided Lie algebra in
the vector space category $_F {\cal M}$
with the ordinary flip  $\tau $ (i.e. $\tau (x\otimes y) = y \otimes x$)
as braiding, although an extension of $L$ may
become a braided Lie algebra in  $_F {\cal M}$.
Furthermore, for the algebra $L = F\{ x \} $  in $(^{F {\bf Z}} {\cal M}, C^r)$
with $q ^2\not=1$ and
$r (k, m) = q ^{km}$ for any $k, m \in {\bf Z}$, the above conclusion holds.

\begin {Definition} \label {1.7}
Let $(L, [\ \  ])$ be a braided  m-Lie algebra in $({\cal C}, C)$.
If $M$ is an object and  there exists a morphism
$\alpha : L \otimes M \rightarrow M$ such that $\alpha ([\ \ ] \otimes M)
= \alpha (L \otimes \alpha ) - \alpha (L \otimes \alpha ) (C \otimes M),$
 then $(M, \alpha )$ is called an $L$-module.
\end {Definition}

\section {Generalized Matrix Braided m-Lie Algebras}
As examples of the braided m-Lie algebras, we  introduce the concepts of
generalized matrix algebras (see \cite {Zh93}) and  generalized matrix
braided m-Lie algebras.  We construct  generalized
classical braided m-Lie algebras $sl_{q, f}( GM_G(A),  F)$ and
$osp_{q, t} (GM_G(A), M, F)$ of generalized matrix algebra $GM_G(A)$.
We show  how generalized matrix  Lie color algebras are related to
Lie superalgebras for any abelian group $G$. That is,
 we establish  the relationship
between generalized matrix  Lie color algebras and Lie superalgebras.

 Let $I$ be a set. For any $ i, j, l, k \in I,$ we choose a vector space $A_{ij}$
 over field $F$ and an  $F$-linear map
$\mu_{ijl}$ from  $A_{ij}\otimes A_{jl}$ into $A_{il}$ (written
$\mu _{ijl} (x, y)=xy)$ such that  $x(yz)=(xy)z$ for any $x\in
A_{ij}$ , $y\in A_{jl} , z\in A_{lk}$.  Let $A$ be the external direct sum
of $\{ A_{ij} \mid i, j\in  I \}$. We define the multiplication in $A$ as
$$xy = \{ \sum _k x_{ik}y_{kj} \}$$
 for any $x=\{x_{ij}\},  y=\{y_{ij}\}\in A$ .
  It is easy to check that  $A$ is an algebra (possibly without unit ).
  We call $A$ a generalized matrix algebra,  or a gm algebra in short,
written as $A=\sum \{A_{ij} \mid i,  j\in  I\}$ or $GM_I(A)$. Every element
in $A$ is called a generalized matrix. We can easily define  gm ideals
and gm subalgebras.  It is easy to define upper triangular generalized
matrices, strictly upper triangular generalized matrices
and diagonal generalized matrices under some total order $ \prec $ of $I$.

\begin {Proposition} \label {2.1}
Let  $A=\sum \{A_{ij} \mid i, j\in I \}$ be  a gm algebra and  $G$
 an abelian group with $G=I$.
 Then $A$ is an algebra graded by $G$ with
$A_g = \sum _{i = j +g} A_{ij}$ for any $g\in G$. In this case,
the gradation is called a generalized matrix gradation, or gm
gradation in short.
\end {Proposition}
{\bf Proof.} For any $g, h \in G$, see that
\begin{eqnarray*}
A_g A_h &=& (\sum _{i = j +g} A_{ij})( \sum _{s = t +h}
A_{st})\\
 &\subseteq & \sum _{i = t+h +g}A_{i,t+h}A_{t+h,t}\\
 &\subseteq & A_{g +h}. \ \ \ \ \ \
 \end{eqnarray*}
Thus $A=\sum\{A_{ij} \mid i, j\in I \}= \sum _{g\in G} A_g$ is a
$G$-graded algebra. $\Box$

If ${\cal  C }$ is a small preadditive category  and
$A_{ij}=Hom_{\cal C}(j, i)$ is a vector space over
$F$ for any $i, j\in  I =$ objects in $ \ {\cal C}$, then we may easily
show that $\sum\{A_{ij} | i, j\in I\}$ is a generalized matrix algebra.
 Furthermore, if $V = \oplus_{g\in G} V_g$ is a
graded vector space over field $F$ with $A_{gh}=Hom_{F}(V_h, V_g)$
for any $g, h \in G$, then we call braided m-Lie  algebra
$A=\sum\{A_{ij} \mid i, j\in G \}$ the general linear braided
m-Lie   algebra, written as $gl (\{V_g\}, F)$. If dim $V_g = n_g <
\infty $ for any $g\in G$, then  $gl (\{V_g\}, F)$ is written as
$gl (\{n_g\}, F)$. Its braided m-Lie  subalgebras are called
linear  braided m-Lie  algebras. In fact, $gl (\{n_g\}, F) = \{ f
\in End _F V \mid ker f \hbox { is finite codimensional } \}.$
When $G$ is finite, we may  view $gl (\{n_g\}, F)$ as a block
matrix algebra over $F.$ When $G=0$, we denote $gl (\{n_g\}, F)$
by $gl (n, F)$, which  is the ordinary ungraded  general linear
Lie algebra.

Assume that $D$ is a directed graph ($D$ is possibly  an infinite
directed graph and also possibly not a simple graph ).  Let $I$
denote the vertex set of $D$, $x_{ij}$ an arrow from $i$ to $j$
and $x=(x_{i_1i_2}, x_{i_2i_3}, \cdots , x_{i_{n-1}i_{n}})$ a path
from $i_1$ to $i_n$  via arrows $x_{i_1i_2},x_{i_2i_3}, \cdots ,
x_{i_{n-1}i_{n}}$. For two paths $x=(x_{i_1i_2},x_{i_2i_3}, \cdots
, x_{i_{n-1}i_{n}})$ and $y=(y_{j_1j_2},y_{j_2j_3}, \cdots ,
y_{j_{m-1}j_{m}})$ of $D$ with $i_n=j_1$, we define the
multiplication of $x$ and $y$ as  $$xy = (x_{i_1i_2}, x_{i_2i_3},
\cdots , x_{i_{n-1}i_{n}}, y_{j_1j_2}, y_{j_2j_3}, \cdots ,
y_{j_{m-1}j_{m}}) .$$ \noindent Let $A_{ij}$ denote the vector
space over field $F$ with basis being all paths from $i$ to $j$,
where $i, j\in I$. Notice that we view every vertex $i$ of $D$ as
a path from $i$ to $i$, written $e_{ii}$ and $e_{ii}x_{ij} =
x_{ij}e_{jj} = x _{ij}$.   We can naturally define a linear  map
from $A_{ij} \otimes A_{jk}$ to $A_{ik}$ as $x\otimes y = xy $ for
any two pathes $x \in A_{ij}, y \in A_{jk}$. We may easily show
that $\sum \{A_{ij}\mid  i, j \in  I \}$ is a generalized matrix
algebra, which is  called a path algebra,  written as $A(D)$ (see,
\cite [Chapter 3]{ARS95}).

\begin {Example} \label {2.3}
There are   finite-dimensional braided  m-Lie  algebras in braided
tensor category $(^{F {\bf Z}_3} {\cal M}, C^r)$ with a primitive
3th root $q$ of 1 and $r (k, m) = q ^{km}$ for any $k, m \in {\bf
Z}_3$. Indeed, for any natural number $n$,  we can construct a
generalized matrix braided  m-Lie algebra $A = \sum \{A_{ij} \mid
i, j \in {\bf Z}_3\}$ such that $dim \ A =n.$

(i) Let $A_{ij} = 0$ when $i\neq j$ but $A_{11}=F$. Thus $dim \ A = 1$.

(ii) Let $A_{ij} = 0$ when $i\neq j$ but $A_{11}=A_{22}=F$. Thus $dim \ A = 2$.

(iii) Let $A_{ij} = 0$ when $i\neq j$ but $A_{11}= A_{22}=A_{33}=F$.
Thus $dim \ A = 3$.

(iv) Let $D$ be a directed graph with vertex set ${\bf Z}_3$ and
only one arrow from 1 to 2. Set $A = A(D).$ It is clear $dim (
A_{ij}) = 0$ when $i\neq j$ but $ dim(A_{11})= dim(A_{22})=dim(A_{33})=dim
( A_{12})=1$. Thus $dim \  A = 4$.

(v) Let $D$ be a directed graph with vertex set ${\bf Z}_3$ and
only two arrows: one  from 1 to 2 and other one from 1 to 3. Set
$A=A(D)$. It is clear  $dim (A_{ij}) = 0 $ but $ dim (A_{11})= dim
( A_{22})=dim ( A_{33})=dim ( A_{12})= dim (A_{13})=1$ Thus $dim \
A = 5$.

vi) Let $D$ be a directed graph with vertex set ${\bf Z}_3$ and
only $n+2$ arrows: one  from 1 to 2,  one from 2 to 3 and the
others from 1 to 3. Set $A= A(D).$ It is clear $dim (A_{ij}) = 0 $
but $ dim (A_{11})= dim ( A_{22})=dim ( A_{33})=dim ( A_{12})= dim
(A_{23})=1$ and $dim (A_{13}) = n +1$.  Thus $dim \ A = n+6 $ for
$n =0, 1, \cdots $.
\end {Example}

Let $A$ be a braided m-Lie algebra, $G$ be an abelian group with a bicharacter
$r $ and $W$ be a vector space over $F$.
\begin{Definition}\label{quantum-trace}
If $f$ is an $F$-linear map from
$GM_G(A)$ to $W$ and satisfies the following

(i)  $f (a) = \sum _{g \in G} r (g, g) f (a_{gg})$

(ii) $f (a_{ij} b_{ji}) =  f (b_{ji} a_{ij})$  for any $a, b \in A$ and
$i, j \in G,$ then $f$ is called a generalized quantum trace function from
gm  algebra $GM_G(A)$ to $W$, written $tr_{q, f}$.
\end{Definition}
Set
$$sl_{q, f}(GM_G(A), F) = \{ a \in GM_G(A)  \mid  tr_{q,f} (a)=0 \}.$$

By computation,
$$tr_{q,f} [a, b] = r(u, u) ^{-1} r(g, g)\sum _{g \in G} (1-  r(u, u)^2 r
(u, g)r (g, u)) tr_{q,f} (b _{g, g +u} a_{g+u, g}) $$
for any homogeneous elements $ a, b \in A, u \in G. $ Thus we get
\begin {Lemma}\label {2.4}
$sl_{q, f}(GM_G(A), F)$ is a braided m-Lie algebra  with
$[A,A] \subseteq sl_{q, f}(GM_G(A), F)$ iff
$$\sum _{g \in G} (1-  r(u, u)^2 r (u, g)r (g, u)) tr_{q, f}
(b _{g, g +u} a_{g+u, g}) =0$$
\noindent for any homogeneous elements $ a, b \in A, u \in G $ with
$\mid $$a$$ \mid = u$ and $\mid $$b$$ \mid = -u.$
If, in addition,  $tr_{q, f} (A_{ij}A_{ji}) \not=0$ for any $i, j \in G$, then
 $sl_{q, f}(GM_G(A), F)$ is a braided m-Lie algebra with
 $[A,A] \subseteq sl_{q, f}(GM_G(A), F)$
 iff $r$ is a skew symmetric bicharacter.
 \end {Lemma}

If $f(a_{gg}) = a_{gg}$ for any $g \in G,$ then from definition
\ref{quantum-trace} $f$ is a generalized quantum
trace  from $GM_G(A)$  to $GM_G(A).$
If $GM_G(A) = M_n(F)$ is the full matrix algebra over $F$ with
$G = {\bf Z}_n$, $r (i, j)=1$ and $f(a_{gg}) = a_{gg} $
for any $i , j , g \in {\bf Z}_n$, then $f$ is the ordinary trace function.
If $GM_G(A) = gl (\{ n_g\}, F)$ and $f (a_{gg})= tr (a_{gg})$
(i.e. $f (a_{gg})$ is the ordinary trace of the matrix $a_{gg}$ )
 for any $g \in G$,
then $f$ is the quantum trace of the graded matrix algebra $gl (\{ n_g\}, F)$.
In this case, $sl_{q, f} (gl (\{n_g\}, F))$ is simply
written as $sl _q (\{n_g\}, F).$

Let $G$ be an abelian group with a  bicharacter
$r $ .  If $t$ is an $F$-linear map from $GM_G(A)$ to $GM_G(A)$
such that  $t (A_{ij}) \subseteq A_{ji}$, $(t(a))_{ij} = t(a _{ji})$ and
$t(ab)= t(b)t(a)$ for any $i, j \in G, a, b \in GM_G(A)$,
 then $t$ is called a generalized transpose on $GM_G(A).$
Given   $ 0\not=M \in GM_G(A) $, for any $u\in G,$
let $ osp _{q, t} (GM_G(A), M, F )_u $ $= \{ a \in GM_G(A) _u \mid
t(a_{u+g, g})M_{u+g,h} = -  r(g,u) M_{g, u+h } a_{u+h, h}
    \hbox { for any } g, h \in G\}$
and   $osp _{q, t}(GM_G(A), M, F) = \oplus  _{u\in G} osp
(GM_G(A), M, F)_u.$

\begin {Lemma}\label {2.5}
$osp_{q, t} (GM_G(A), M, F)$ is a braided m-Lie  algebra iff
$$\sum _{g \in G} (1-  r (u, v)r (v, u)) M_{g, h +u +v} a _{h +u +v, v+h}
b _{v +h, h} =0$$
\noindent for any homogeneous elements $ a, b \in A, u, v \in G $ with
$\mid $$a$$ \mid = u$ and $\mid $$b$$ \mid = v.$
If, in addition,  for any $i, j,  k\in G$,  $a_{ij} \not=0$ implies
$a_{ij}A_{jk} \not=0$, then $osp_{q, t} (GM_G(A), M, F))$ is a braided m-Lie
algebra  iff $r$ is a skew symmetric bicharacter.
\end {Lemma}
{\bf Proof.}  Obviously,  $osp_{q, t} (GM_G(A), M, F)_u$
is a subspace.  It remains to check that $osp _{q, t} (GM_G(A), M, F)$ is
closed under bracket operation. For any
$a \in  osp_{q, t} (GM_G(A), M, F)_u$, $b \in
osp _{q, t} (GM_G(A), M, F)_v$ and  $ u, v, g, h \in G,$ set $w = u+v.$ See that
\begin {eqnarray*}
t([a, b] _{w+g, g}) M _ {w+g, h}
&=& t((ab - r (v, u)ba) _{w + g, g})M _ {w+g, h} \\
&=&t(b_{v + g, g})t(a_ {w+g, g+v})M _ {w+g, h} \\
 &{\ \ \ }& -
r (v, u)t(a_{u + g, g})t(b_ {w+g, g+u})M _ {w+g, h}\\
&=& r (g +v, u)r(g, v)M _{g, h +w}b _{h+w, u+h}a _{u+h, h} \\
&{\ \ \ }&- r(v,u)r (g +u, v)r(g, u)M _{g, h +w}a _{h+w, v+h}b _{v+h, h}, \\
-  r(g,w) M_{g, w+h } [ a, b]_{w+h, h}
 &=& - r(g, w) M_{g, w +h} a _{w+h, v+h}b _{v+h, h}\\
& &- r(g, w) r (v , u) M _{g, w+h} b _{w+h, u +h} a _{u +h, h}.
\end {eqnarray*}
Thus $[a, b] \in osp _{q, t} (GM_G(A), M, F)_w $ iff
$(r (u, v)r(v, u) -1) M_{g, w+h}
 a_{w+h, v+h} b _{v+h, h} =0$ for any $g,  h \in G.$
 $\Box$

We now  consider   $sl_{q, f}( GM_G(A),  F)$  and
$osp_{q, t} (GM_G(A), M, F)$ when the  bicharacter $r$ is skew symmetric.
In this case, they become  Lie color algebras in
$(^{FG} {\cal M}, C^r)$,  called special gm Lie color algebra and
ortho-symplectic gm Lie color algebra, respectively.

It is well-known that a Lie color algebra $(^{FG}{\cal M}, C^r)$
with finitely generated $G$ is related to a Lie super algebra by
\cite [Theorem 2] {Sc79}.  We now show  how the above gm Lie color
algebras  are related to  Lie superalgebras for any abelian group
$G$.

Let $G$ be an abelian group with a skew symmetric bicharacter $r $. Set
$G_ {\bar 0} = \{g \in G \mid r (g, g )=1 \}$ and
$G_ {\bar 1} = \{g \in G \mid r (g, g )=-1 \}$. We define a new bicharacter:
$r_0 (g, h) = -1$ for $g, h\in G_{\bar 1}$ and $r_0 (g, h) =1$ otherwise.
It is clear that $r_0$ is a bicharacter too.
 Obviously, $(r_0)_0 = r_0$ for any skew symmetric bicharacter $r$ on $G$.
 For convenience,
let $L^s $ denote the Lie superalgebra $ L= L_{\bar 0} \oplus
L_{\bar 1}$ for Lie color algebra $L$ in $(^{FG}{\cal M}, C^r)$,
 where $L_{\bar 0} = \oplus _{ i \in G_{\bar 0}} L_i$ and
 $L_{\bar 1} = \oplus _{ i \in G_{\bar 1}} L_i$ with
 $[ x, y ] = xy - r_0(\mid $$y$$ \mid  , \mid $$x$$ \mid )yx $ for any
$x \in L_g, y\in L_h, g, h \in G$ (see \cite [Page 718] {Sc79}).

Let $A = \sum \{A_{ij} \mid i, j \in G\} = GM_G(A)$ be a
generalized matrix algebra. Set $B_{\bar i , \bar j} = \sum _{g
\in G_{\bar i}, h \in G_{\bar j}} A_{gh}$ for any $\bar i , \bar j
\in {\bf Z}_2$ and $B = \sum \{B_{\bar i , \bar j } \mid  \bar i ,
\bar j \in {\bf Z}_2\}$. We denote  the generalized matrix algebra
$GM_{{\bf Z}_2}(B)$ by $ GM_{{\bf Z}_2}(A^s) $.  For $a\in
GM_G(A),$  if $b_{\bar i , \bar j} = \sum _{g \in G_{\bar i}, h\in
G_{\bar j} } a_{gh}$ for any $\bar i , \bar j \in {\bf Z}_2$, then
we denote the element $b$ by $a^s$. When $G= {\bf Z}_2$ with $r
(g,h) = (-1) ^{gh}$ for any $g, h \in {\bf Z}_2,$ we denote
$tr_{q, f}$ by  $ tr_{s, f}$ and $osp _{q, t}$ by  $ osp _{s, t}.$
We have
\begin {Theorem}\label {2.6}

(i)  $sl_{q, f}( GM_G(A),  F)^s  = sl _{s, f} (GM_{{\bf Z}_2}(A^s), F).$

(ii) $osp_{q,t}(GM_G(A), M, F) ^s = osp _{s, t} (GM_{{\bf Z}_2}(A^s), M^s, F).$
 \end {Theorem}
{\bf Proof.} (i) For any $a \in  GM_G(A)$, see that
 \begin{eqnarray*}
tr _{q, f} (a) &=& \sum _{g \in G_{\bar 0}} r_0 (g, g) tr_{q,f} (a_{gg}) +
 \sum _{g \in G_{\bar 1}} r_0 (g, g) tr_{q,f} (a_{gg})\\
&=& tr_{s, f} (a).
\end {eqnarray*}
This completes the proof of (i).

(ii) For any $a \in (osp _{q, t} (GM_G(A), M, F)^s )_{\bar 0}$ with
$a \in osp _{q, t} (GM_G(A), M, F)_u$ and
$u \in G_{\bar 0},$ let $\bar i , \bar j \in {\bf Z}_{2}$ and  we see that
\begin{eqnarray*}
t(a _{ \bar i , \bar  i}) M_{\bar  i , \bar  j}
&=& \sum_{g, h \in G_ {\bar i}, k \in G_{\bar j}} t(a _{gh}) M _{gk} \\
&=& - \sum_{ h \in G_ {\bar i}, k \in G_{\bar j}} r_0
(h, u) M _{h, k +u} a_{k+u, k} \\
&=& - M_{\bar i, \bar j } a_{\bar j , \bar j }.
\end{eqnarray*}
\noindent This shows $a \in osp _{s, t} (GM_{{\bf Z}_2}(A^s), M^s, F)$.
We can similarly prove the others.  $\Box$

In fact, the relations \ \ \ (i) and (ii)\ \ \  above define a
$G$-grading of Lie superalgebras \ \ \ \ \ \ \ \ \ \ $sl _{s, f}
(GM_{{\bf Z}_2}(A^s), F)$ and $ osp _{s, t} (GM_{{\bf Z}_2}(A^s),
M^s, F)$, respectively.

We may apply the above results to $gl (\{n _g\}, F)$ with the ordinary
quantum trace $tr_q (a) = \sum _{g\in G} r(g,g) tr(a_{gg})$ and the
ordinary transpose  $t (a) = a'.$ In this case,
$sl _{q, f}$ and $osp _{q, f}$ are denoted by $sl _q $ and $osp _q,$
respectively. We have
\begin {Corollary}\label {2.7}
(i)  $(sl_{q}( gl (\{V_g\}, F))^s = sl _{s} (V_{\bar 0}, V_{\bar 1}, F).$

(ii) $ (osp_q gl(\{V_g\}, F))^s = osp _s (V _ {\bar 0 }, V_{\bar 1},  M^s, F).$
 \end {Corollary}

\section {Braided m-Lie Algebras in the Yetter-Drinfeld Category and
Their Representations } In this section, we  give another class of
braided m-Lie algebras. We shall show that representations of an
algebra $A$ associated to a braided m-Lie algebra $L$ are also
representations of $L$. Furthermore, we show that if $(M, \psi)$ is a
faithful representation of $L$, then
 representations of $End _F M$ are also representations of  $L$.

The category is the Yetter-Drinfeld category $(^B_B{\cal YD}, C)$,
where $B$ is a finite dimensional Hopf algebra and $C$ is
a braiding with $C(x,y) = \sum (x_{(-1)}\cdot y)\otimes x_{(0)}$ for any
$x\in M, y\in N$.

We use the Sweedler's notation for  coproducts and comodules, i.e.
$$\Delta (a) = \sum _a a_1 \otimes a_2  \hbox { \ \ \ and  \ \ \ \ }
\psi (x) = \sum _x x_{(-1)} \otimes x_{(0)}$$ when $a\in H$ a coalgebra and
$x\in M$ a left $H$-comodule.

\begin {Lemma} \label {3.1}
(see \cite [Lemma 2.1 and Lemma 2.3 (iv)]{ZZH03})

(i) If $(V, \alpha_V, \phi _V)$ and $(W, \alpha_W, \phi _W)$ are two
Yetter-Drinfeld modules over $B$ with dim $B<\infty $, then
 $Hom _F (V, W)$ is a Yetter-Drinfeld module under the following module
 operation and comodule operation:
 $(b \cdot f)(x) = \sum b_1\cdot f(S(b_2)\cdot x )$ and
$\phi (f) = (S^{-1} \otimes \hat \alpha ) (b_B \otimes f )$,
 where  $ \hat \alpha $ is defined by
$(b^* \cdot f)(x) = <b^*, x_{(-1)} S(f (x_{(0)})_{(-1)})>_{ev}
(f (x_{(0)}))_{(0)}$ for any $x\in V, f \in Hom _F (V,W), b^* \in B^*.$
Here $b_B$ denotes a coevoluation and $<, >_{ev}$ an evoluation of $B.$

(ii) If  $(M, \alpha_M, \phi _M)$ is  a Yetter-Drinfeld modules over $B$,
$End _F M$ is an algebra in $(^B_B {\cal YD}, C)$
 \end {Lemma}

{\bf Proof.}
 (i)  It is clear  that
$\sum f _{(-1)} f _{(0)}(x) = \sum (f (x_{(0)}))_{(-1)} S^{-1}(x_{(-1)})
\otimes (f (x_{(0)}))_{(0)} $ for any $x\in V, f \in Hom _F (V, W), b \in B.$
Using this, we can show that $Hom _F (V, W)$ is a $B$-comodule.
Similarly, we can show that $Hom _F (V, W)$ is a $B$-module. We now show
\begin {eqnarray} \label {e3.1}
\sum (b \cdot f )_{(-1)} \otimes (b \cdot f )_{(0)}
= \sum b_1 f _{(-1)}S(b_3) \otimes b_2 \cdot f _{(0)}
\end {eqnarray}
for any $f \in Hom _F(V, W), b \in B$. For any $x\in V,$  see that
\begin {eqnarray*}
 \sum (b \cdot f )_{(-1)} \otimes (b \cdot f )_{(0)}(x)&=&
 \sum b_1 (f (S(b_4))\cdot x_{(0)})_{(-1)} S(b_3)S^{-1}(x_{(-1)})\\
&{\ }& \otimes b_2 \cdot (f (S(b_4))\cdot x_{(0})_{(0)}, \\
b_1 f _{(-1)}S(b_3) \otimes (b_2 \cdot f _{(0)})(x)&=& b_1 f_{(-1)} S(b_4)
\otimes  b_2 \cdot f_{(0)} ((S(b_3)\cdot x))\\
&=& \sum b _1 (f (S(b_4) \cdot x _{(0)}))_{(-1)} S(b_3)S^{-1}(x_{(-1)})
b_5S(b_6) \\
&{\ }& \otimes b_2 \cdot (f (S(b_4) \cdot x _{(0)}))_{(0)} \\
&=& \sum b_1 (f (S(b_4))\cdot x_{(0)})_{(-1)} S(b_3)S^{-1}(x_{(-1)})\\
&{\ }& \otimes b_2 \cdot (f (S(b_4))\cdot x_{(0)})_{(0)}.
\end {eqnarray*}
Thus (\ref {e3.1}) holds and $Hom _F (V, W)$ is a Yetter-Drinfeld module.

(ii) Let $E = End _F M$ and $m$ denote the multiplication of $E$.
 Now we show that $m$ is a homomorphism of  $B$-comodules.
It is sufficient to show
\begin {eqnarray} \label {e3.11}
\sum _{fg}(fg )_{(-1)} \otimes (fg )_{(0)}
= \sum _{f, g} f_{(-1)} g_{(-1)} \otimes f_{(0)}g_{(0)}.
\end {eqnarray}
for any $f, g \in E$. Indeed, for any $x\in M,$  see that
 \begin {eqnarray*}
\sum _{fg}(fg )_{(-1)} \otimes (fg )_{(0)}(x)  &=&
\sum _{x } (fg (x_{0}))_{(-1)} S^{-1}(x_{(-1)})
\otimes (fg (x_{0}))_{(0)},\\
\sum _{f, g} f_{(-1)} g_{(-1)} \otimes f_{(0)}g_{(0)} (x)
 &=& \sum _{f, x} f_{(-1)} (g (x_{(0)}))_{(-1)} S^{-1} (x_{(-1)})
 \otimes f_{(0)} ((g (x_{(0)}))_{(0)}) \\
 &=& \sum _{f, x} (f((g (x_{(0)}))_{(0)}))_{(-1)} \\
&{\ }& S^{-1}((g (x_{(0)}))_{(0)}))_{(-1)}{}_2)
(g (x_{(0)}))_{(-1)} {}_1 S^{-1} (x_{(-1)}) \\
&{\ }& \otimes  (f((g (x_{(0)}))_{(0)}))_{(0)} \\
&=&   (fg(x_{(0)}))_{(-1)} S^{-1}(x_{(-1)})\otimes (fg(x_{(0)}))_{(0)}.
\end {eqnarray*}
Thus (\ref {e3.11}) holds. Similarly, we can show that $m$ is a homomorphism
of  $B$-modules.
 $\Box$

\begin {Example} \label {3.2}
Let  $(M, \alpha_M, \phi _M)$ be  a Yetter-Drinfeld modules over $B$ and
$L$ a subobject
of  $E = End _F M$. If $L$ is  closed under operation  $[\ \ ] = m -
m C_{L, L},$ then $(L, [\  \ ] )$ is a braided  m-Lie  subalgebra of $E^-$.
\end {Example}

By Lemma \ref {3.1}, we may define representations of braided m-Lie algebras
in the Yetter-Drinfeld category.

\begin {Definition} \label {3.3}
Let $(L, [\ \  ])$ be a braided  m-Lie algebra in $(^B_B{\cal YD}, C)$.
If $M$ is an object in $(^B_B{\cal YD}, C)$ with  morphism
$\psi : L \rightarrow End _F M$
such that $\psi $ is a homomorphism of braided m-Lie algebras, i.e.
 $\psi [\  \ ] = m (\psi \otimes \psi ) -
m (\psi \otimes \psi )  C_{L, L},$
then $(M, \psi)$ is called a representation of $(L, [\ \  ])$.
\end {Definition}
Obviously, $(M, \psi )$ is a  representation of $L$  iff  $(M,
\alpha )$ is an $L$-module (see definition \ref {1.7}), where the
relation between two operations  is $\alpha (a,  x)  = \psi
(a)(x)$ for any $a\in L, x\in M$.

\begin {Proposition} \label {3.4}
Let $(L, [\ \  ])$ be an object in  $(^B_B{\cal YD}, C)$
with a morphism $[\ \ ] : L \otimes L \rightarrow L$.
If $M$ is an object in $(^B_B{\cal YD}, C)$ with monomorphism
 $\psi : L \rightarrow End _F M$
such that $\psi [\  \ ] = m (\psi \otimes \psi ) -
m (\psi \otimes \psi )  C_{L, L},$
then  $L$ is a braided m-Lie algebra and $(M, \psi)$ is a faithful
representation of $(L, [\ \  ])$.
\end {Proposition}
{\bf Proof.} By Lemma \ref {3.1} $E= End _F M$ is an algebra in
$(^B_B {\cal YD}, C)$. By Definition \ref {1.1},
$(L, [\ \ ])$ is a braided m-Lie algebra. $\Box$

\begin {Proposition} \label {3.5}
 Let $(L, [ \ \ ])$ be  a braided m-Lie algebra in $(_B^B {\cal YD}, C)$
 induced by multiplication of $A$ through $\phi .$

(i) If $( M , \psi )$  is a representation of algebra $A$, then
$(M, \psi \phi)$ is a representation of $L.$

(ii) If $( M , \psi )$  is a representation of braided m-Lie
algebra $A{}^-$, then $(M, \psi \phi)$ is a representation of $L.$

(iii) If   $( M , \psi )$  is  a faithful representation of $L$
and  $(N, \varphi )$ is a representation of algebra $End_F M$,
then
 $(N, \varphi \psi  )$   is a representation of  $L$.
\end {Proposition}

{\bf Proof.} (i) and (ii) follow from  Definition \ref {3.3}.

 (iii) By Lemma \ref {3.1}, $E = End_F M$ is
an algebra in $(^B_B{\cal YD}, C)$. Thus $L$ is a braided m-Lie
algebra induced by $E$ through $\psi$. By Proposition \ref {3.5},
we complete the proof. $\Box$

Let $(L, [\ \ ])$ be  a braided  m-Lie subalgebra of the path algebra
$(F(D, \rho ))^-$ with relations,  then a  representation
$(V, f )$ of $D$ with $f_\sigma =0$ for any $\sigma \in \rho$ is also a
representation of $L$
(see \cite [Proposition II.1.7]{ARS95}).

\vskip 0.5cm {\bf Acknowledgement }: The authors are financially
supported by Australian Research Council. S.C.Z thanks the
Department of Mathematics, University of Queensland for
hospitality.

\end {document}